\documentclass[12pt]{article}
\usepackage{amsmath,amssymb,amsfonts,amsthm,xcolor}
\usepackage[margin=1in]{geometry}
\usepackage{subcaption}
  
\usepackage{graphicx}

\theoremstyle{plain}

\title{Revisiting Cases 2 and 11 of the Map Color Theorem}
\author{Timothy Sun\\Department of Computer Science\\San Francisco State University}
\date{}

\newcommand{\Z}{\mathbb{Z}}

\begin{document}

\maketitle

\begin{abstract}
In 1968, Ringel and Youngs solved the remaining cases of the orientable Map Color Theorem by finding genus embeddings of the complete graphs $K_n$, for sufficiently large $n \equiv 2, 8, 11 \pmod{12}$. Following the approach previously explored by the author for $n \equiv 8 \pmod{12}$, we aim to streamline their constructions for $n \equiv 2, 11 \pmod{12}$ by finding families of current graphs with simpler patterns for the arc labelings. 
\end{abstract}

\section{Introduction}

The orientable Map Color Theorem, which establishes the chromatic number of every closed orientable surface besides the sphere, boils down to finding minimum genus embeddings of the complete graphs $K_n$. The proof is split depending on the residue $n$ modulo 12, about which Mohar and Thomassen \cite{MoharThomassen} wrote, ``for the most complicated [residues], no short proofs are known.''

As chronicled in Ringel and Youngs \cite{RingelYoungs}, outside of a few small cases, the three residues $n \equiv 2, 8, 11 \pmod{12}$ were the last to be solved, and these constructions are collected in Section 7 of Ringel \cite{Ringel-MapColor}. For $k = 0, \dotsc, 11$, the problem of finding a genus embedding of $K_n$ for each $n \equiv k \pmod{12}$ is referred to as \emph{Case $k$}. Our primary focus is on Case 11, which is arguably the most complicated residue at present.

The only known proofs of the orientable Map Color Theorem rely on the theory of \emph{current graphs}, which are arc-labeled, embedded graphs that produce symmetric embeddings of covering graphs. Current graphs can be used to generate triangular embeddings, but for some residues $n$ modulo 12, the complete graphs $K_n$ do not triangulate any orientable surface. The standard approach uses current graphs to generate triangular embeddings of near-complete graphs, after which the missing edges are added using handles and other local operations. The latter process is referred to as the \emph{additional adjacency} part of the construction.

In the proof presented in Section 7 of Ringel \cite{Ringel-MapColor}, the additional adjacency steps require specific edge flips, which induce constraints on the requisite current graphs. Ringel and Youngs' strategy for finding such current graphs starts by fixing the underlying graph, which contains a long ``ladder'' subgraph, and labeling the edges outside of the ladder according to the constraints. Afterwards, the remaining currents are assigned to the ladder, often with the help of some other labeling problem (e.g. graceful labelings of path graphs). 

Jungerman (for example, in \cite{Jungerman-KnK2}) approached this problem in reverse order. First, a repeating pattern of currents, derived from the simplest graceful labeling of the path graph, is assigned to the ladder. The length of the ladder and the labels can be chosen so that the problem of finding an \emph{infinite} family of current graphs is reduced to a \emph{finite} problem of routing and assigning the remaining, constant-sized set of edges and currents. In the author's experience, Jungerman's method has been more fruitful in finding novel families of current graphs. Ringel and Youngs's solution to Case 8 was previously simplified by the author \cite{Sun-FaceDist} using this technique, and we will demonstrate similar success with Cases 2 and 11. Though Ringel and Youngs felt ``that no possible further improvement will ever be made on Case 11'' \cite[p.viii]{Ringel-MapColor}, certain aspects of our new solution are simpler. 

The work presented here represents possibly the final chapter in the author's attempts at simplifying the proof of the orientable Map Color Theorem. Previously, the author found improvements to Cases 0 \cite{Sun-K12s}, 1 \cite{Sun-FaceDist}, 6 \cite{Sun-Minimum}, and 8 \cite{Sun-FaceDist}. The current graphs we present here are of index 1, i.e., ones where the embedding has one face, but we note that there are also higher-index approaches to Cases 2 \cite{Jungerman-KnK2, Sun-Index2}, 8, and 11 \cite{Sun-Minimum}. 

\section{Graph embeddings}

For background on topological graph theory, see Gross and Tucker \cite{GrossTucker} and Mohar and Thomassen \cite{MoharThomassen}. 

In this work, all graph embeddings are orientable and cellular. Let $S_k$ denote the orientable surface of genus $k$, the $k$-holed torus. Given a graph $G = (V,E)$ embedded in an orientable surface $S_k$, the connected components of $S \setminus \phi(G)$ are called \emph{faces} and are denoted by the set $F$. Cellular embeddings satisfy the \emph{Euler polyhedral equation}
$$|V| - |E| + |F| = 2-2k.$$
The \emph{(orientable) genus} $\gamma(G)$ is the smallest value $k$ such that $G$ has an embedding in $S_k$. If $G$ is simple, connected, and has at least three vertices, the Euler polyhedral equation implies a lower bound on its genus:
$$\gamma(G) \geq \left\lceil \frac{|E|-3|V|+6}{6} \right\rceil.$$
The genus of the complete graphs $K_n$, due to Ringel, Youngs, and others \cite{Ringel-MapColor}, matches this lower bound:
$$\gamma(K_n) = \left \lceil \frac{(n-3)(n-4)}{12} \right\rceil,$$
for $n \geq 3$. 

For each edge $e \in E$, we orient the edge arbitrarily to obtain two arcs $e^+$ and $e^-$ that point in opposite directions. The set of all such arcs is denoted by $E^+$. An embedding in an oriented surface can be described by a \emph{rotation system}, where each vertex is assigned a \emph{rotation}, a cyclic permutation of all of the arcs leaving the vertex. We can obtain a rotation system from such an embedding by considering the clockwise ordering of all arcs leaving the vertex. One can recover the cellular embedding from a rotation system by face-tracing, and the genus of the surface can be calculated using the Euler polyhedral equation. When the graph is simple, a rotation can also be specified as a cyclic permutation of the neighbors of the vertex. 

\section{Current graphs}

In this work, a \emph{current graph} is a graph embedded in an orientable surface whose arcs are labeled with elements of an abelian group $\Gamma$, known as the \emph{current group}. The arc-labeling $\alpha\colon E^+ \to \Gamma$ satisfies $\alpha(e^+) = -\alpha(e^-)$ for each edge $e \in E$.  The \emph{excess} of a vertex is defined to be the sum of all the currents on the arcs entering the vertex, and if the excess is 0, we say that \emph{Kirchhoff's current law} is satisfied at that vertex. Given a face-boundary walk $(e_1^\pm, e_2^\pm, \dotsc)$, the \emph{log} of the walk $(\alpha(e_1^\pm), \alpha(e_2^\pm), \dotsc)$ replaces each arc with its label. 
 
The \emph{index} of a current graph refers to the number of faces in its embedding. We follow the approach described by Ringel and Youngs \cite{RingelYoungs-Case2, RingelYoungs-Case11}, where the current graphs are of index 1 and the current group is $\Z_{12s+6}$. Most of the vertices satisfy Kirchhoff's current law, and those that do not are called \emph{vortices} (with one exception, as explained below). Each face corner incident with a vortex is labeled by a letter. Ringel and Youngs used the following kinds of vortices:
\begin{itemize}
\item[(V1)] Vertex of degree 1 with excess generating $\Z_{12s+6}$.
\item[(V2)] Vertex of degree 1 with excess generating the index 2 subgroup.
\item[(V3)] Vertex of degree 3 with excess generating the index 3 subgroup, and the incoming currents are all congruent to $j \pmod{3}$, where $j \not\equiv 0 \pmod{3}$.
\end{itemize}
With these vortices in mind, the current graphs of Ringel and Youngs's method satisfy a standard set of properties:
\begin{itemize}
\item[(C1)] Each nonzero element of $\Z_{12s+6}$ appears in the log exactly once.
\item[(C2)] Vortices are of type (V1), (V2), or (V3). 
\item[(C3)] The order 2 element is on an arc incident with an unlabeled vertex of degree 1. 
\item[(C4)] All other vertices are unlabeled, have degree 3, and satisfy Kirchhoff's current law.
\end{itemize}

By property (C3), the order 2 element appears twice consecutively in the log, which will induce faces of length 2. We suppress these faces, replacing them by a single edge, and we follow the conventions of recording the order 2 element only once in the log, not considering the degree-1 endpoint as a vortex, and omitting that endpoint from our drawings. Vortices of type (V1) and (V3) generate one and three Hamiltonian $(12s+6)$-sided faces, respectively. Vortices of type (V2) each induce two $(6s+3)$-sided faces, where one face is incident with the even elements of $\Z_{12s+6}$, and the other is incident with the odd elements. We also follow the convention where the letters are incorporated into the log at their respective face corners in the face-boundary walk. 

The derived embedding is initially of a complete graph with vertex set $\Z_{12s+6}$. The rotation at vertex $i \in \Z_{12s+6}$ is generated by temporarily ignoring the letters in the log and adding $i$ to each element in the log. For each of the Hamiltonian faces, we subdivide them with the corresponding letters so that the rotation at vertex 0 is identical to the log. For the two $(6s+3)$-sided faces generated by a vortex of type (V2), if the vortex label is, say, $y$, then the two subdivision vertices are called $y_0$ and $y_1$, where $y_i$ is adjacent to vertex $i$, for $i = 0,1$. After subdividing long faces and suppressing two-sided faces, we obtain a triangular embedding. 

Vertices which are elements of the current group are called \emph{numbered} vertices, and those that come from vortex labels are called \emph{lettered} vertices. The numbered vertices are already pairwise adjacent, by property (C1), so the additional adjacency step in the constructions aims to identify the pairs of vertices from (V2) vortices and add the missing edges between all of the lettered vertices.

\begin{figure}[ht]
\centering
\includegraphics[scale=1]{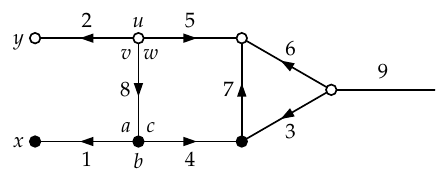}
\caption{A current graph with current group $\Z_{18}$.}
\label{fig-case2-s1}
\end{figure}

Since we will later need to examine the rotation systems of the derived embeddings in detail, we give a small example in Figure \ref{fig-case2-s1}, where black and white vertices represent clockwise and counterclockwise rotations, respectively. It contains each type of vortex and represents the smallest case of the infinite family we will later see in Figure \ref{fig-case2}. The log of its face-boundary walk is
$$\arraycolsep=4pt\begin{array}{rrrrrrrrrrrrrrrrrrrrrrrrrrrrrrrrrrrrrrrr}
(9 & 6 & 13 & u & 2 & y & 16 & v & 8 & c & 4 & 7 & 12 & 3 & 14 & b & 1 & x & 17 & a & 10 & w & 5 & 11 & 15).
\end{array}$$

Consequently, the rotations at the first few numbered vertices are
$$\arraycolsep=4pt\begin{array}{clcccccccccccccccccccccccr}
0. & (9 & 6 & 13 & u & 2 & y_0 & 16 & v & 8 & c & 4 & 7 & 12 & 3 & 14 & b & 1 & x & 17 & a & 10 & w & 5 & 11 & 15) \\
1. & (10 & 7 & 14 & v & 3 & y_1 & 17 & w & 9 & a & 5 & 8 & 13 & 4 & 15 & c & 2 & x & 0 & b & 11 & u & 6 & 12 & 16) \\
2. & (11 & 8 & 15 & w & 4 & y_0 & 0 & u & 10 & b & 6 & 9 & 14 & 5 & 16 & a & 3 & x & 1 & c & 12 & v & 7 & 13 & 17) \\
3. & (12 & 9 & 16 & u & 5 & y_1 & 1 & v & 11 & c & 7 & 10 & 15 & 6 & 17 & b & 4 & x & 2 & a & 13 & w & 8 & 14 & 0) \\
4. & (13 & 10 & 17 & v & 6 & y_0 & 2 & w & 12 & a & 8 & 11 & 16 & 7 & 0 & c & 5 & x & 3 & b & 14 & u & 9 & 15 & 1) \\
\vdots
\end{array}$$

As we read from top to bottom, the positions corresponding to a (V3) vortex cycle through its three vertex labels. The order of the labels depends on the rotations of the vortices and whether the incoming currents are congruent to $1$ or $2 \pmod{3}$. 

\begin{figure}[ht]
\centering
\includegraphics[scale=1]{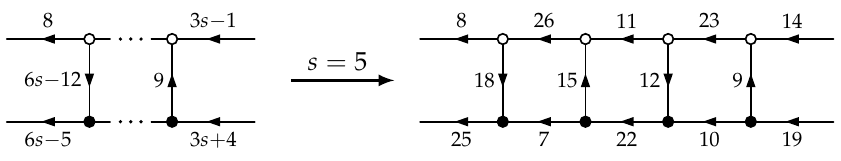}
\caption{A ladder defined for all $s \geq 3$ and its specification for $s = 5$.}
\label{fig-ladder}
\end{figure}

Our infinite families of current graphs contain a ladder fragment like the one shown in Figure \ref{fig-ladder}. The ladder's ``rungs'' alternate in direction, and the currents are consecutive multiples of 3. As mentioned earlier, this assignment is based off of the simplest graceful labeling of a path graph, and is simpler than any of the ones described in Section 7 of Ringel \cite{Ringel-MapColor}.  

\section{Case 11}

Ringel and Youngs \cite{RingelYoungs-Case11} solved Case 11 using current graphs containing two vortices of type (V1), labeled $x$ and $y$, and a vortex of type (V3), whose corners are labeled $a$, $b$, and $c$. The derived embedding is a triangular embedding of the graph $K_{12s+11}-K_5$. For expository purposes, we attempt to express their additional adjacency step as generally as possible and reinterpret their two handle operations as a single modification that increases the genus by 2.

Given an edge $e$ incident with two triangular faces, removing $e$ and replacing it with the other ``diagonal'' $e'$ of the resulting quadrangular face is called an \emph{edge flip}. If $e'$ connects two vertices that are already adjacent, we keep the graph simple by deleting that other edge. This triggers a sequence of edge flips that we can engineer to connect two specific nonadjacent vertices. 

\begin{figure}[ht]
\centering
\includegraphics[scale=1]{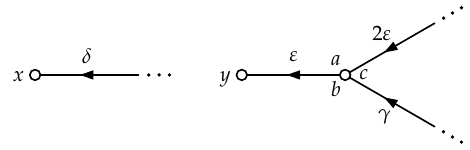}
\caption{The initial set of constraints on vortices for Case 11.}
\label{fig-case11boundary}
\end{figure}

\begin{figure}[ht]
\centering
    \begin{subfigure}[b]{0.99\textwidth}
        \centering
        \includegraphics[scale=1]{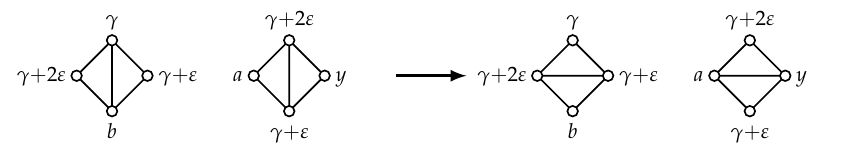}
        \caption{}
    \end{subfigure}
    \begin{subfigure}[b]{0.99\textwidth}
        \centering
        \includegraphics[scale=1]{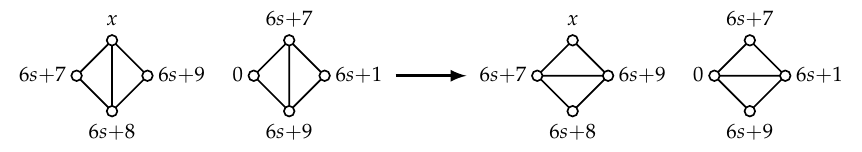}
        \caption{}
    \end{subfigure}
\caption{A sequence of edge flips that are possible on any current graph satisfying the initial constraints (a), and a sequence specific to our new family of current graphs (b).}
\label{fig-case11flips}
\end{figure}

We begin by constraining the vortices, as illustrated in Figure \ref{fig-case11boundary}. Let $\delta$ and $\varepsilon$ be the excesses of the (V1) vortices $x$ and $y$, respectively. We require that $\varepsilon \equiv 2 \pmod{3}$, that the (V3) vortex is adjacent to vortex $y$, and that the log is of the form
$$\begin{array}{rrrrrrrrrrrrrrrrrrrrrrrrrrrrrrrrrrrrrrrr}
(\dotsc & 2\varepsilon & a & \varepsilon & y & -\varepsilon & b & \dotsc & c & \dotsc & x & \dotsc)
\end{array}$$
Let $\gamma$ be the current of the remaining arc entering the (V3) vortex. Since this current must satisfy $\gamma \equiv -\varepsilon \equiv 1 \pmod{3}$, the sequence of edge flips starting by deleting the edge $(\gamma, b)$, as seen in Figure \ref{fig-case11flips}(a), ends by adding the edge $(a, y)$. We also flip the edge $(0, -\varepsilon)$ to obtain $(b,y)$. 

\begin{figure}[ht]
\centering
    \begin{subfigure}[b]{0.35\textwidth}
        \centering
        \includegraphics[scale=1]{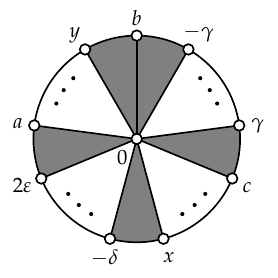}
        \caption{}
    \end{subfigure}
    \begin{subfigure}[b]{0.35\textwidth}
        \centering
        \includegraphics[scale=1]{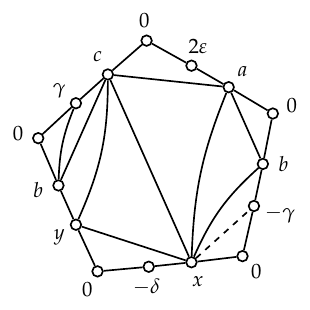}
        \caption{}
    \end{subfigure}
\caption{The handle operation merges the five shaded triangles (a) into a 15-sided face (b) where all the missing edges can be incorporated.}
\label{fig-case11mod}
\end{figure}

The rotation at vertex $0$ is now of the form
$$\begin{array}{rrrrrrrrrrrrrrrrrrrrrrrrrrrrrrrrrrrrrrrr}
(b & -\gamma & \dotsc & \gamma & c & \dotsc & x & -\delta & \dotsc & 2\varepsilon & a & \dotsc & y)
\end{array}$$
We consider the following five subsequences: $b$, $-\gamma \dotsc \gamma$, $c \dotsc x$, $-\delta \dotsc 2\varepsilon$, and $a \dotsc y$. By rearranging these subsequences so that the rotation becomes
$$\begin{array}{rrrrrrrrrrrrrrrrrrrrrrrrrrrrrrrrrrrrrrrr}
(b & a & \dotsc & y & -\delta & \dotsc & 2\varepsilon & c & \dotsc & x & -\gamma & \dotsc & \gamma),
\end{array}$$
the five shaded faces in Figure \ref{fig-case11mod}(a) ``in between'' the subsequences are merged into the 15-sided face
$$[0, c, \gamma, 0, b, y, 0, -\delta, x, 0, -\gamma, b, 0, a, 2\varepsilon].$$

\begin{figure}[ht]
\centering
    \begin{subfigure}[b]{0.99\textwidth}
        \centering
        \includegraphics[scale=1]{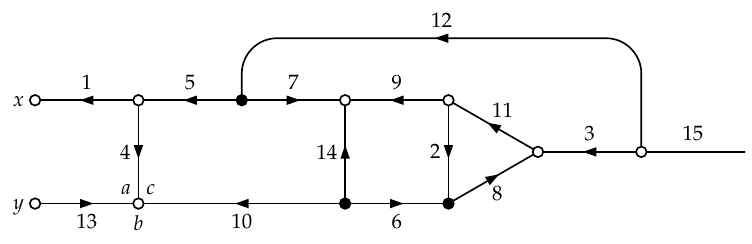}
        \caption{}
    \end{subfigure}
    \begin{subfigure}[b]{0.99\textwidth}
        \centering
        \includegraphics[scale=1]{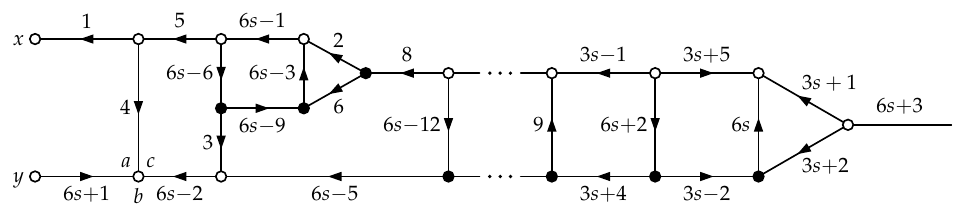}
        \caption{}
    \end{subfigure}
\caption{A family of current graphs with current group $\Z_{12s+6}$, defined for (a) $s = 2$ and (b) $s \geq 3$.}
\label{fig-case11}
\end{figure}

Following the solid lines in Figure \ref{fig-case11mod}(b), all of the missing edges can be drawn inside of this face except $(0, -\varepsilon)$. There are three quadrangular faces, $[0, c, a, 2\varepsilon]$, $[0, -\delta, x, y]$, and $[0, -\gamma, b, x]$, where adding one of the diagonals initiates a sequence of edge flips that can potentially end with adding the final missing edge. The current graphs of Ringel and Youngs \cite{RingelYoungs-Case11} allow for such a sequence starting with the edge $(2\varepsilon, c)$. Our family of current graphs in Figure \ref{fig-case11}, which has parameters $\gamma = 6s{-}2$, $\delta = 1$, and $\varepsilon = 6s{+}5$, uses the edge $(-\gamma, x)$ (the dashed line in Figure \ref{fig-case11mod}(b)), instead. The resulting sequence of edge flips is shown in Figure \ref{fig-case11flips}(b), where the quadrangle in the second edge flip has the opposite orientation when $s = 2$. 

We note that these current graphs are very similar to those used for Case 8 in Sun \cite{Sun-FaceDist}. In particular, the ladder with rungs $6s-12, \dotsc, 9$ is nearly identical, with only the directions of the arcs changed, and there is a sequence of edge flips involving a (V1) vortex of excess $\pm 1$ and an edge with currents $\pm 2$. Whether our family is simpler than that of Ringel and Youngs \cite{RingelYoungs-Case11} is a matter of taste. On one hand, our infinite family is only defined for $s \geq 3$, while Ringel and Youngs's family starts at $s = 2$, though with slight differences depending on the parity of $s$. On the other hand, the rungs on our ladder form a single arithmetic sequence (as opposed to two separate sequences) and the vertices in the ladder satisfy Kirchhoff's current law over $\Z$ (rather than just over $\Z_{12s+6}$), unifying the construction with other known index 1 solutions. Furthermore, our sequence of edge flips is shorter by one step, and the orientations of the edge flips are consistent for all $s \geq 3$. 

\section{Case 2}

\begin{figure}[ht]
\centering
\includegraphics[scale=1]{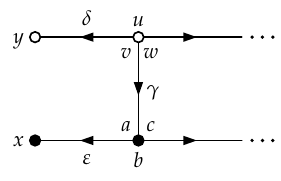}
\caption{Vortices for Case 2, where $x$ and $y$ are of type (V1) and (V2), respectively.}
\label{fig-case2-boundary}
\end{figure}

The arc-labeling in Ringel and Youngs's solution for Case 2 \cite{RingelYoungs-Case2} reuses much of their labeling in Case 11 that we circumvented in the previous section. Thus, it seems natural to simplify Case 2, as well. Their initial conditions on the current graphs, shown in Figure \ref{fig-case2-boundary}, add eight additional vertices (the two vertices $y_0$ and $y_1$ that result from the (V2) vortex are ultimately identified). No additional constraints are needed. 


\begin{figure}[ht]
\centering
    \begin{subfigure}[b]{0.35\textwidth}
        \centering
        \includegraphics[scale=1]{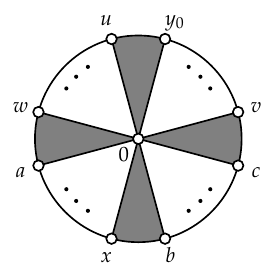}
        \caption{}
    \end{subfigure}
    \begin{subfigure}[b]{0.45\textwidth}
        \centering
        \includegraphics[scale=1]{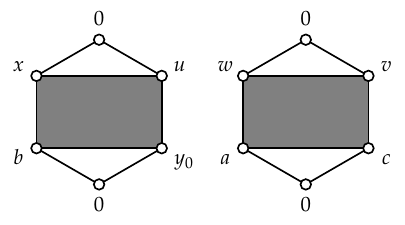}
        \caption{}
    \end{subfigure}
\caption{The handle operation merges pairs of shaded triangles (a) into two 6-sided faces (b).}
\label{fig-case2mod}
\end{figure}

Since the additional adjacency step is almost identical to that in Ringel and Youngs \cite{RingelYoungs-Case2}, we only sketch the details. The rotation at vertex $0$ is of the form
$$\arraycolsep=4pt\begin{array}{rrrrrrrrrrrrrrrrrrrrrrrrrrrrrrrrrrrrrrrrrrr}
(\dotsc & u & \delta & y_0 & \dotsc & v & \gamma & c & \dotsc & b & \varepsilon & x & \dotsc & a & -\gamma & w & \dotsc)
\end{array}$$
First, the edges $(0, \delta)$, $(0, \gamma)$, $(0, \varepsilon)$, and $(0, -\gamma)$ are flipped and replaced with the edges $(u, y_0)$, $(v, c)$, $(b, x)$, and $(a, w)$, respectively. Then, similar to Case 11, we group the remaining neighbors of vertex $0$ into the subsequences $y_0 \dotsc v$, $c \dotsc b$, $x \dotsc a$, and $w \dotsc u$. We change the rotation at 0 by permuting the subsequences:
$$\arraycolsep=4pt\begin{array}{rrrrrrrrrrrrrrrrrrrrrrrrrrrrrrrrrrrrrrrrrrr}
(y_0 & \dotsc & v & w & \dotsc & u & x & \dotsc & a & c & \dotsc b).
\end{array}$$
Now, as shown in Figure \ref{fig-case2mod}, there are two six-sided faces $[0, x, b, 0, y_0, u]$ and $[0, w, a, 0, c, v]$. After adding the edges $(x,u)$, $(b, y_0)$, $(w, v)$, and $(a, c)$ inside of these faces, following Section 7.5 of Ringel \cite{Ringel-MapColor}, one may connect the two shaded faces $[x, b, y_0, u]$ and $[w, a, c, v]$ with three handles to add all of the remaining missing edges between these eight vertices. In Sun \cite[Fig.\ 9]{Sun-Minimum}, this process is interpreted as gluing a certain genus embedding of $K_8$ to these two quadrangular faces. After all of these modifications, the rotations at numbered vertices besides $0$ are the same as in the original derived embedding, except possibly with $0$ missing. 

\begin{figure}[ht]
\centering
\includegraphics[scale=1]{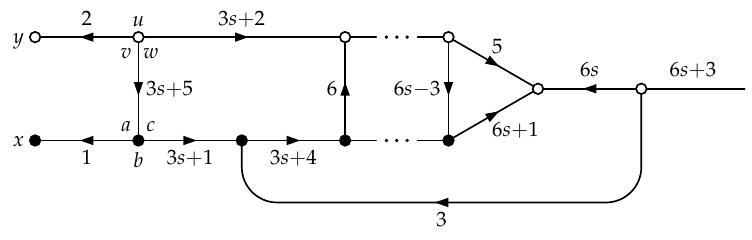}
\caption{A family of current graphs with current group $\Z_{12s+6}$, defined for all $s \geq 1$.}
\label{fig-case2}
\end{figure}

At this point, we need to restore the edges $(0, \delta)$, $(0, \gamma)$, $(0, -\gamma)$, and $(0, \varepsilon)$ and identify $y_0$ with $y_1$. The element $\varepsilon$ is guaranteed to be odd, since it is the excess of a vortex of type (V1). Of the remaining three vertices $\delta$, $-\gamma$, and $\gamma$, there is always a triangular face containing two of those vertices, where the third vertex, which we call $v$, is adjacent to $0$. If not, then each triple from $\{\delta, -\gamma, \gamma, \varepsilon\}$ would form a triangular face. Since the embedding is in a manifold, this would imply that each of these vertices has degree 3, a contradiction. For concreteness, our current graphs in Figure \ref{fig-case2} have the parameters $\gamma = 3s+5$, $\delta = 2$, and $\varepsilon = 1$, and there is a triangle of the form $[v, -\gamma, \delta]$, where $v = c$ for $s = 1$, and $v = 9s+10$ for $s \geq 2$.

Assume, without loss of generality, that the rotation at vertex $v$ is of the form 
$$\arraycolsep=4pt\begin{array}{rrrrrrrrrrrrrrrrrrrrrrrrrrrrrrrrrrrrrrrrrrr}
(\gamma & \dotsc & \delta & -\gamma & \dotsc & 0 & u_1 & \dotsc & u_i),
\end{array}$$
where $u_1$ and $u_i$, $i \geq 1$, are any other vertices. Like with the previous handle operations, we permute the subsequences $\gamma \dotsc \delta$, $-\gamma \dotsc 0$, and $u_1 \dotsc u_i$ to obtain
$$\arraycolsep=4pt\begin{array}{rrrrrrrrrrrrrrrrrrrrrrrrrrrrrrrrrrrrrrrrrrr}
(-\gamma & \dotsc & 0 & \gamma & \dotsc & \delta & u_1 & \dotsc & u_i).
\end{array}$$

\begin{figure}[ht]
\centering
    \begin{subfigure}[b]{0.35\textwidth}
        \centering
        \includegraphics[scale=1]{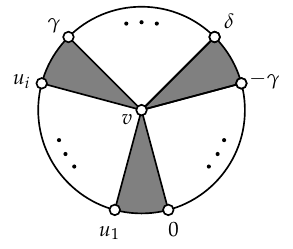}
        \caption{}
    \end{subfigure}
    \begin{subfigure}[b]{0.35\textwidth}
        \centering
        \includegraphics[scale=1]{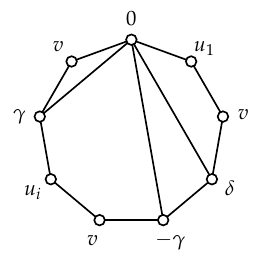}
        \caption{}
    \end{subfigure}
\caption{The handle operation merges the three shaded triangles (a) into one (b).}
\label{fig-p3mod}
\end{figure}

The three shaded faces in Figure \ref{fig-p3mod}(a) are merged into one face incident with $0$, $\delta$, $\gamma$, and $-\gamma$. Then, the edges $(0, \delta)$, $(0, \gamma)$, and $(0, -\gamma)$ are drawn inside of this face. 

Finally, since there are edges $(0, y_0)$ and $(\varepsilon, y_1)$, we can attach a handle near these edges to add the final missing edges $(0, \varepsilon)$ and $(y_0, y_1)$. We contract the edge $(y_0, y_1)$ to obtain an embedding of the complete graph $K_{12s+14}$. As calculated by Ringel \cite[p.116]{Ringel-MapColor}, this is a minimum genus embedding, for all $s \geq 1$. 

Once again, the arc-labelings are simpler than that of Ringel and Youngs \cite{RingelYoungs-Case2}, but there is another minor disagreement between the smallest case and the general case, namely the identity of vertex $v$. We leave open the possibility of further unifying the solutions to both of these residues. 

\bibliographystyle{alpha}
\bibliography{biblio}

\begin{thebibliography}{RY69b}

\bibitem[GT87]{GrossTucker}
Jonathan~L. Gross and Thomas~W. Tucker.
\newblock {\em {Topological Graph Theory}}.
\newblock Wiley \& Sons, 1987.

\bibitem[Jun75]{Jungerman-KnK2}
Mark Jungerman.
\newblock The genus of {$K_n-K_2$}.
\newblock {\em Journal of Combinatorial Theory, Series B}, 18(1):53--58, 1975.

\bibitem[MT01]{MoharThomassen}
Bojan Mohar and Carsten Thomassen.
\newblock {\em Graphs on {S}urfaces}, volume~10.
\newblock Johns Hopkins University Press, 2001.

\bibitem[Rin74]{Ringel-MapColor}
Gerhard Ringel.
\newblock {\em {Map Color Theorem}}, volume 209.
\newblock Springer Science \& Business Media, 1974.

\bibitem[RY68]{RingelYoungs}
Gerhard Ringel and J.W.T. Youngs.
\newblock {Solution of the Heawood map-coloring problem}.
\newblock {\em Proceedings of the National Academy of Sciences},
  60(2):438--445, 1968.

\bibitem[RY69a]{RingelYoungs-Case11}
Gerhard Ringel and J.W.T. Youngs.
\newblock Solution of the {H}eawood map-coloring problem --- {C}ase 11.
\newblock {\em Journal of Combinatorial Theory}, 7(1):71--93, 1969.

\bibitem[RY69b]{RingelYoungs-Case2}
Gerhard Ringel and J.W.T. Youngs.
\newblock {Solution of the {H}eawood map-coloring problem --- {C}ase 2}.
\newblock {\em Journal of Combinatorial Theory}, 7(4):342--352, 1969.

\bibitem[Sun19]{Sun-K12s}
Timothy Sun.
\newblock A simple construction for orientable triangular embeddings of the
  complete graphs on $12s$ vertices.
\newblock {\em Discrete Mathematics}, 342(4):1147--1151, 2019.

\bibitem[Sun20]{Sun-Minimum}
Timothy Sun.
\newblock Simultaneous current graph constructions for minimum triangulations
  and complete graph embeddings.
\newblock {\em Ars Mathematica Contemporanea}, 18(2):309--337, 2020.

\bibitem[Sun21]{Sun-FaceDist}
Timothy Sun.
\newblock Face distributions of embeddings of complete graphs.
\newblock {\em Journal of Graph Theory}, 97(2):281--304, 2021.

\bibitem[Sun24]{Sun-Index2}
Timothy Sun.
\newblock Jungerman ladders and index 2 constructions for genus embeddings of
  dense regular graphs.
\newblock {\em European Journal of Combinatorics}, 120:103974, 2024.

\end{thebibliography}

\end{document}